\documentclass{article}
\usepackage{amsmath,amsfonts, amssymb,verbatim}
\newtheorem{Theo}{Theorem}  
\newtheorem{pkt}{}[subsection]  
\newcommand{\bpk}{\begin{pkt}\rm }  
\newcommand{\epk}{\end{pkt}}  
\newtheorem{Lemm}{Lemma}[subsection]  
   
\newtheorem{Prop}{Proposition} 
   
\newtheorem{Cor}{Corollary}   
   
\newcommand{\pf}{{\bf Proof}\ }   
   
\newcommand{\HH}{{\cal H}}

 \newcommand{\W}{{\bf W}}

 \newcommand{\al}{\alpha}

\newcommand{\M}{{\bf M}}   
\newcommand{\A}{{\cal A}}

\newcommand{\p}{{\bf p}} 
\newcommand{\R}{{\mathbb R}}   
   
\newcommand{\Z}{{\mathbb Z}}   
\newcommand{\N}{{\mathbb N}}  
\newcommand{\C}{{\mathbb C}}  
  
\newcommand{\ee}{\end{equation}}  
\newcommand{\be}{\begin{equation}}

\newcommand{\pr}{{\rm pr} }

\newcommand{\ssn}{\section}

\newcommand{\bl}{\begin{Lemm}}   
   
\newcommand{\el}{\end{Lemm}}   
   
\newcommand{\bt}{\begin{Theo}}   
   
\newcommand{\et}{\end{Theo}}   
   
\newcommand{\bp}{\begin{Prop}}   
   
\newcommand{\ep}{\end{Prop}}   
   
\newcommand{\inv}{^{-1}}   
   
\newcommand{\bc}{\begin{Cor}}   
   
\newcommand{\ec}{\end{Cor}}   
  
\newcommand{\lb}{\label}  
\newcommand{\sgn}{{\rm ang}}  
\newcommand{\bd}{{\rm bd}}

\newcommand{\ra}{\rangle}  
\newcommand{\la}{\langle}  
   
\newcommand{\subs}{\subseteq}

\newcommand{\qed}{$\Box$ \\ \\}   
\newcommand{\F}{{\rm F}}  
\newcommand{\T}{{\rm T}}  
\newcommand{\FF}{{\bf U}} 
\newcommand{\GGG}{{\bf V}}

\newcommand{\eps}{\epsilon}   
\newcommand{\Ac}{{\cal A}} 
\newcommand{\XX}{{\frak X}}

\newcommand{\xx}{{\bf x}} 
\newcommand{\yy}{{\bf y}} 
\newcommand{\zz}{{\bf z}} 
 
\newcommand{\Y}{{\bf Y}} \newcommand{\X}{{\bf X}} 
\newcommand{\ZZ}{{\bf Z}}    
\newcommand{\PPP}{{\rm P}} 
 
 \newcommand{\DDD}{{\bf D}}

\newcommand{\f}{{\bf u}} \newcommand{\g}{{\bf v}}  
\newcommand{\ww}{{\bf w}}

\begin{document}
\title{Non-commutative Zariski geometries and their classical limit}
\author{Boris Zilber\\ University of Oxford}
\date{5 July 2007}
\maketitle
\abstract{We undertake a case study of two series of 
nonclassical Zariski geometries. We show that 
these geometries can be realised as representations of 
certain noncommutative $C^*$-algebras and introduce a natural limit construction which for each of the 
two series produces a classical $U(1)$-gauge field over a $2$-dimensional Riemann surface. }

\section{Introduction}
The notion of a Zariski geometry was introduced in [HZ] as a model-theoretic
generalisation of objects of algebraic geometry and compact complex manifolds.

The main result of paper [HZ] was the 
classification of non linear (non locally modular) irreducible Zariski 
geometries of dimension one. The initial hope that every such geometry
is definably isomorphic to an algebraic curve over an algebraically
closed field $\F$  had to be corrected in the course of the study.
The final classification theorem states that given a non linear 
irreducible  1-dimensional Zariski geometry $M$ there is an
algebraically closed field $\F$ definable in $M$ and an algebraic
curve $C$ over $\F$ such that $M$ is a  finite cover of $C(\F),$ that
is there is a Zariski continuous map  $\p: M\to C(\F)$
which is a surjection with finite fibres. 

The paper [HZ] also provides a class of examples that demonstrates
that in general we can not hope to reduce $\p$ to a bijection. 
 Given a smooth algebraic curve $C$ with a big enough group $G$ of regular
automorphisms with a nonsplitting finite extension $\tilde G,$
one can produce a ``smooth'' irreducible Zariski curve $\tilde C$
along with a finite cover $\p: \tilde C\to C$ and $\tilde G$ its group
of Zariski-definable automorphisms. 
 
Typically $\tilde C$ can not be identified with any algebraic curve 
because $\tilde G$ 
is not embeddable into the group
of regular automorphisms of an algebraic curve ([HZ], section 10).

Taking into account known reductions of covers we can say that the
above construction of $\tilde C$ is essentially the only way to
produce a nonclassical Zariski curve. In other words, a general
Zariski curve essentially looks like $\tilde C$ above.

A simple example of an unusual group $\tilde G$
for such a $\tilde C,$ used in [HZ], is the class-$2$-nilpotent group of two generators
$\f$ and $\g$ with the central commutator $[\f,\g]$ of finite
order $N.$ The correspondent $G$ is then the free abelian group on two
generators. 
One can identify this $\tilde G$ as the quotient of the integer
Heisenberg group $H_3(\Z)$ by the subgroup of its centre of index $N.$

Also, since the group of regular isomorphisms of the smooth curve $C$ must be
infinite, we have very little freedom in choosing
$C;$
it has to be either the affine line over $\F,$ or the torus $\F^*,$ or
an elliptic curve.

This paper undertakes the case study of the geometries of the 
corresponding $\tilde C$ for $C$ an algebraic torus and an affine line. \\

The most comprehensive modern notion of a {\em geometry} is based on
 the consideration of a {\em coordinate algebra} of the geometric object.
The classical meaning of a coordinate algebra comes from
the algebra of  {\em coordinate functions} on the  object,
that is, in our case, functions $\psi: \tilde C(\F)\to \F$
of a certain class. The most natural algebra of functions 
seems to be the algebra $\F[\tilde C]$ of Zariski continuous (definable) 
functions. 
But by the virtue of the construction $\F[\tilde C]$ is naturally
isomorphic to the 
 $\F[C],$ the algebra of regular functions on the algebraic curve $C,$
that is the only geometry which we see by looking into
$\F[\tilde C]$ is the geometry of the algebraic curve $C.$
To see the rest of the structure we had to extend $\F[\tilde C]$ by 
introducing {\em semi-definable} functions, which satisfy certain {\em
equations} but are not uniquely defined by these equations. The
$\F$-algebra of $\HH(\tilde C)$ of semi-definable functions contains
the necessary information  about $\tilde C$ but is not canonically
defined. On the other hand it is possible to define an $\F$-algebra 
$\Ac(\tilde C)$ of linear operators on   $\HH(\tilde C)$ in a
canonical way, depending on $\tilde C$ only. 
We proceed with this construction for both
examples and write down explicit lists of generators and defining
relations for algebras $\Ac(\tilde C).$ One particular type of
a semi-definable function
which we call {\em  $*$-functions}, of a clearly non-algebraic nature, plays a special role. The $*$-function induces 
 an involution $*$ on $\Ac.$
We show, for  $\F=\C,$ that $\Ac$ thus gets the structure of 
a $C^*$-algebra, that is the involution $*$ associates with any $X\in \Ac$
its formal {\em adjoint} operator $X^*$ satisfying usual formal requirements.
 Moreover there is an $\Ac$-submodule 
$\HH_+$ of  $\HH(\tilde C)$ with
an inner product for which $*$ does indeed define adjoint operators.

Our first main theorem states  that there is a reverse 
canonical
 construction which  recovers $\tilde C$ from  the algebra $\Ac$ uniquely.     
The points of $\tilde C$ correspond to one-dimensional eigenspaces (states)
of certain self-adjoint operators, relations on $C$ correspond to ideals
 of cartesian powers of a commutative subalgebra of $\Ac$ and operations $\f$ and $\g$ correspond
 naturally to actions of certain
 operators of $\Ac$ on the states.
This scheme is strikingly similar to the
operator representations of quantum mechanics.  Note that this construction
is similar but not identical with the one we used in [Z1].

The final section of the paper concentrates on understanding the limit
of the structures $\tilde C=\tilde C_N,$ depending on $N$ by the  construction
of $\tilde G,$ as $N$
tends to infinity. Amomg many possible ways to define the notion of
the limit we found metric considerations most relevant. It turns out
possible, when $\F=\C,$  to consider metric on each 
$\tilde C_N$ and to use correspondingly the notion of
Hausdorff limit. Our main result in this section states that,
for both types of examples, the Hausdorff limit $\tilde
C_{\infty}$ of  $\tilde C_N,$ as $N$ tends to infinity, is the structure identified as the
principal $U(1)$-bundle over a Riemann surface   with $\f$ and
$\g$ defining a connection (covariant derivative) on the bundle. In
physicists' terminology this is a gauge field 
 with a nonzero curvature 
 (see e.g. [DFN] or [S]).

Combining with the results of the previous section one could speculate
that $\tilde C_N$ are quantum deformations of the classical
structure on $\tilde C_{\infty},$ and conversely, the latter is the classical
limit of the quantum structures. \\ \\

I am very grateful to Yang-Hui He and Mario Serna for helping me to clarify some of physicists' terminology that I am using in the paper. 

\section{  Non-algebraic Zariski geometries}\lb{znorv}  
\bpk \lb{tilde}
\bt    
 There exists an  irreducible  pre-smooth Zariski structure   
(in particular of dimension $1$) which is not interpretable in an algebraically closed field.\et  
  
{\bf The construction}\\

 Let $\M=(M,{\cal C})$ be an irreducible  pre-smooth Zariski structure,\\  
$G\le {\rm ZAut}\, M$ (Zariski-continuous bijections) acting freely on $M$ and for some $\tilde G$ with {\bf finite $H:$}  
$$1\to H\to \tilde G\to^{\pr} G\to 1.$$  
  
Consider a set $X\subs M$  of representatives of $G$-orbits:\ \ for each   
$a\in M,$\ \ $G\cdot a\,\cap X$ is a singleton.  
  
Consider the formal set $$\tilde M(\tilde G)=\tilde M=\tilde G\times X$$  
and the projection map  
$$\p: (g,x)\mapsto \pr(g)\cdot x.$$  
Consider also, for each $f\in \tilde G$ the function  
$$f: (g,x)\mapsto (fg,x).$$  
We thus have obtained  
the structure $$\tilde \M=(\tilde M,\ \{ f\}_{f\in \tilde G}\cup   
\p\inv({\cal C}))$$  
on the set $\tilde M$ with relations induced from $\M$ together with  
maps $\{ f\}_{f\in \tilde G}.$ We set the closed subsets of $\tilde M^n$ to be exactly those which are definable by positive quantifier-free formulas with   
parameters.  Obviously, the structure $\M$ and the map  
$\p: \tilde M\to M$ are definable in $\tilde \M.$ Since, for each $f\in \tilde G,$  
$$\forall v\ \p f(v)=f\p(v)$$  
the image $\p(S)$ of a closed subset  $S\subs \tilde M^n$ is closed in $\M.$  
We define $\dim S:= \dim \p(S).$\\  
   
{\bf Lemma 1} { \em   
 The isomorphism type of $\tilde\M$ is determined by $M$ and $\tilde G$ only.  
The theory of $\tilde \M$ has quantifier elimination.  
 $\tilde \M$  
is an irreducible pre-smooth Zariski structure.}     
  
 Proof. One can use obvious   
automorphisms of the structure to prove quantifier elimination. The statement of the claim then follows by checking the definitions. The detailed proof is given in [HZ] Proposition 10.1.\\

{\bf Lemma 2}. {\em Suppose $H$ does not split, that is for every proper $G_0<\tilde G$   
  $$G_0\cdot H\neq \tilde G.$$  
Then, every equidimensional Zariski expansion $\tilde \M'$ of $\tilde \M$  
is irreducible.}\\  
  
\pf Let $C=\tilde M'$ be an   
$|H|$-cover of the variety $M,$ so $\dim C=\dim M$ and $C$ has at  
most $|H|$ distinct irreducible components, say   
$C_i,$ $1\le i\le n.$   
For generic $y\in M$ the fiber $\p\inv(y)$ intersects every $C_i$ (otherwise   
$\p\inv(M)$ is not equal to $C$).   
  
 Hence $H$ acts transitively on the set of irreducible components.  
So, $\tilde G$ acts transitively on  the set of irreducible components, so   
   the setwise stabiliser $G^0$ of $C_1$  in $\tilde G$ is of index $n$  
in $\tilde G$ and also $H\cap\tilde G^0$ is of index $n$ in $H.$   
Hence, $$\tilde G=G^0\cdot H,\mbox{ with }H\nsubseteq G^0$$  
 contradicting our assumptions. \qed

{\bf Lemma 3}. {\em $\tilde G\le {\rm ZAut}\, \tilde \M,$ that is  
$\tilde G$ is a subgroup of the group  of 
Zariski-continuous  bijecions of $\tilde \M.$}  
\\  
  
\pf Immediate by  construction.\qed

{\bf Lemma 4}  {\em Suppose $\M$ is a rational or elliptic  curve (over an algebraically closed field $\F$ of characteristic zero), $H$ does not split, $\tilde G$ is nilpotent and for some big enough integer $\mu$  
there is a non-abelian subgroup $G_0\le \tilde G$ $$|\tilde G:G_0|\ge \mu.$$  
Then  $\tilde \M$ is not interpretable in an algebraically closed field. }\\  
  
\pf First we show.  
  
Claim. Without loss of generality we may assume that $\tilde G$ is infinite.  
  
Recall that $G$ is a subgroup of the group ${\rm ZAut}\, M$ of rational (Zariski) automorphisms of $M.$ Every algebraic curve is birationally equivalent to a smooth one,  
so $G$ embeds into  the group of  
birational transformations of a smooth rational curve or an elliptic curve. Now remember that any birational  
transformation of a smooth algebraic curve is biregular.  
If $M$ is rational then  the group ${\rm ZAut}\, M$  
is ${\rm PGL}(2,{\F}).$ Choose a semisimple (diagonal)   
$s\in {\rm PGL}(2,{\F})$   
 be an automorphism of infinite order such that $\la s\ra \cap G=1$ and $G$ commutes with $s.$ Then we can replace $G$ by $G'=\la G, s\ra$ and $\tilde G$ by $\tilde G'=\la\tilde G, s\ra$ with the trivial action of $s$ on $H$. One can easily see from the construction that the $\tilde M'$ corresponding to $\tilde G'$ is the same as $\tilde M,$ except for the new definable bijection corresponding to $s.$  
  
 We can use the same argument when $M$ is an elliptic curve, in which case   
 the group of automorphisms of the curve   is given as a semidirect product  
 of a finitely generated abelian group   
 (complex multiplication) acting on the group on the elliptic curve  $E(\F).$\\  
  
Now, assuming that $\tilde M$ is definable in an algebraically closed field  
$\F'$ we will have that $\F$ is definable in $\F'.$   
It is known to imply that $\F'$ is definably isomorphic to $\F,$ so   
we may assume that $\F'=\F.$   
  
Also, since $\dim \tilde M=\dim M=1,$ it follows that $\tilde M$ up to finitely many points is in a bijective definable  
correspondence with a smooth algebraic curve, say $C=C(\F).$  
  
$\tilde G$ then by the argument above is embedded into the group of rational automorphisms of $C.$

   The automorphism group is finite if genus of the curve is 2 or higher,   
so by the Claim    we can have only rational or elliptic curve for $C.$

Consider first the case when $C$ is rational.  
   The automorphism  group then   
is ${\rm PGL}(2,{\F}).$ Since $\tilde G$ is nilpotent  
its Zariski closure in ${\rm PGL}(2,{\F})$ is an infinite   
nilpotent group $U.$ Let $U^0$ be the connected component of $U,$ which is     
a normal subgroup of finite index. By Malcev's Theorem (see [Merzliakov], 45.1) there is a number $\mu$ (dependent only on the size of the matrix group in question but not on $U$)  
such that some normal subgroup $U^0$ of $U$ of index at most $\mu$   
 is a subgroup of the  unipotent group    
$$\left( \begin{array}{ll}   
                                1\ \ \ z\\  
                                0\ \ \ 1   
                          \end{array}\right)  
$$   
this is Abelian, contradicting the assumption that $\tilde G$ has no abelian subgroups of index less than $\mu.$   
  
In case $C$ is an elliptic curve  the group of automorphisms   
is  a semidirect product  
 of a finitely generated abelian group   
 (complex multiplication) acting freely on the abelian group of the elliptic curve.   
This group has no nilpotent non-abelian subgroups. This finishes the proof of   
Lemma 4 and of the theorem.\qed

In general it is harder to analyse the situation when $\dim M>1$ since the  
group of birational automorphisms is not so immediately reducible to the group of biregular automorphisms of a smooth variety in higher dimensions. But    
nevertheless the same method can prove the useful fact that the construction produces examples essentially of non algebro-geometric nature.\\    
  
\bp  (i) Suppose  $M$ is an abelian  
 variety, $H$ does not split and $\tilde G$  is nilpotent not abelian. Then $\tilde M$ can not be   
an algebraic variety with $p:\ \tilde M\to M$ a regular map.

(ii) Suppose $M$ is the (semi-abelian) variety $(\F^{\times})^n.$ Suppose also that $\tilde G$ is nilpotent and for some big enough integer $\mu=\mu(n)$  
has no  abelian subgroup $G_0$ of index bigger than $\mu.$ Then  
 $\tilde M$ can not be an  
 algebraic  variety with $p:\ \tilde M\to M$ a regular map.\ep  
\pf (i) If $M$ is an abelian variety 
and $\tilde M$ were algebraic,  
the map $p:\tilde M\to M$ has to be unramified since all its fibers are of the same order (equal to $|H|$). Hence $\tilde M$  
being a finite unramified cover must have the same unversal cover as $M$ has. So, $\tilde M$ must be  an abelian variety as well.    
The group of automorphisms of an abelian variety $\A$ without complex multiplication is the abelian group   
  $\A(\F).$ The contradiction.  
  
(ii) Same argument as in (i) proves that $\tilde M$ has to be isomorphic to $(\F^{\times})^n.$ The Malcev theorem cited above finishes the proof.\qed

\bp Suppose $M$ is an $\F$-variety and, in the  
construction of $\tilde M,$ the group $\tilde G$ is finite.  Then $\tilde M$ is definable in any expansion of the field $\F$ by a total linear order.   
  
In particular, if $M$ is a complex variety, $\tilde M$ is definable in the reals.\ep  
\pf Extend the ordering of $\F$ to a linear order of $M$ and define $$S:=\{ s\in M:\ s=\min\, G\cdot s\}.$$   
The rest of the construction of $\tilde M$ is definable.\qed  
  
{\bf Remark} In other known examples of non-algebraic $\tilde M$ (with  
$G$ infinite) $\tilde M$ is still definable in any expansion of the  
field $\F$ by a total linear order. In particular, for the example  
considered in this paper, see section~\ref{coord}.  
\\  
  
{\bf Problem} (i) Classify Zariski structures definable in the reals.  
  
(ii)  Classify Zariski structures definable in the reals as a smooth real manifold.   
  
(iii) Find new Zariski structures definable in $\R_{an}$ as a smooth real manifold.   
  
\epk

\ssn{Examples}\lb{s2.1}

 Let $N$ be a positive   
integer and  $\F$  an 
algebraically closed field  of characteristic prime to $N.$ 
 Consider the groups given by generators and defining relations,   
$$G=\la u,v:\ uv=vu \ra,$$ $$\tilde G=\tilde G_N=\la \f,\g:\  
[\f,[\f,\g]]=[\g,[\f,\g]]=1=[\f,\g]^N\ra,$$  
where $[\f,\g]$ stands for the commutator $\g\f\g\inv \f\inv.$ 
 
We will consider two examples of the construction of a 
 one-dimensional 
 $\tilde M$ from an algebraic curve $M$ 
using the groups $G$ and $\tilde G.$ By section \ref{znorv} $G$ is 
 going to be a subgroup of the group of rational automorphisms of $M,$ 
 so $M$ has to be of genus $0$ or $1.$ In our examples $M$ is the algebraic torus 
$\F^*$ and 
the affine line $\F.$

\subsection{The $N$-cover of the affine line. } 

\bpk \lb{intrT2}  
We assume here that the characteristic of $\F$ is $0.$ 
 
Let 
$a,b\in \F$ be additively independent.

$G$ acts on $\F:$  
$$u x=a+x,\ v x=b+x.$$  
  
Taking $M$ to be $\F$ this determines,  
by subsection~\ref{znorv}, a presmooth  non-algebraic Zariski curve  
$\tilde M$ which from now on we denote $\PPP_N,$ and $P_N$ will stand
for the universe of this structure.\\

The correspondent definition for the covering map $p: \tilde M\to M=\F$ then gives us  
\be \lb{p12} \p(\f t)=a+\p(t),\ \ \p(\g t)=b+\p(t).\ee  
  
\epk

\bpk \lb{local2} {\bf Semi-definable functions on $\PPP_N.$}\\

{\bf Lemma} {\em  
There are functions $y$ and $z$   
$${P_N}\to \F   $$ 
satisfying  the following {\bf functional equations}, 
 for any $t\in P_N,$  
\be \lb{y}\yy^N(t)=1,\ \yy(\f t)=\eps \yy(t),\  \yy(\g t)=\yy(t)\ee  
\be \lb{z} \zz^N(t)=1,\ \zz({\f t})=\zz(t),\ \zz(\g t)=\yy(t)\inv\cdot \zz(t).\ee  
} 
 
\pf  Choose a subset $S\subs M=\F$ of representatives of $G$-orbits, 
that is $\F=G+S.$ By the construction in section~\ref{znorv}  
we can identify $P_N=\tilde M$ with 
$\tilde G\times S$ in such a way that $\p(g, s)=\pr(g)+s.$ 
This means that,  
for any $s\in S,$ a $t$ in $\tilde G\cdot s$ is of the form 
$t=\f^m\g^n[\f,\g]^{\ell}\cdot s$ and 
  
$$\p(\f^m\g^n[\f,\g]^{\ell}\cdot s):= 
ma+nb+ s.$$ 
Set   
$$\yy(\f^m\g^n[\f,\g]^{\ell}\cdot s):=\eps^{m}$$ 
$$\zz(\f^m\g^n[\f,\g]^{\ell}\cdot s):=\eps^{l}.$$ 
 
This satisfies (\ref{y}) and (\ref{z}). 
\qed  
 
{\bf Remark.} Notice, that it follows from (\ref{p12})-(\ref{z}): 
\begin{enumerate} 
\item $\p$ is surjective and $N$-to-$1,$ with fibres of the form 
$$\p\inv(\lambda)= H t,\ \  
H=\{ [\f,\g]^{\ell}: 0\le l<N\}.$$ 
 
\item $\yy([\f,\g] t)=\yy(t),$  
 
\item $\zz([\f,\g] t)=\eps \zz(t).$  
 
\end{enumerate} 
 
\epk 
 
\bpk \lb{angdef2}  
 
Denote $\F[N]=\{ \xi\in \F: \xi^N=1\}$ and define the {\bf band function} on $\F $  
  as a function 
$\bd: \F  \to \F[N].$ 
  
Set for $\lambda\in \F$  
$$\bd(\lambda)= y(t),\mbox{ if }\p(t)=\lambda,$$ 
This is well-defined by the remark in \ref{local2}.

  
Acting by $\f$ on $t$ and using  (\ref{p12}) and (\ref{y}) we have 
\be \lb{sgn10} \bd(a+\lambda)=\eps \bd\, \lambda.\ee  
Acting by $\g$ we obtain 
\be \lb{sgn20}\bd( b+\lambda)=\bd\, \lambda.\ee 

{\bf Remark} In a more general context we are going to call the
band function and the angular function of the next section
{\bf $*$-functions}, explaining the reasons for this at the end of this 
section.\\

{\bf Proposition} {\em  The structure $\PPP_N$ is definable in  
$$(\F,+,\cdot, \bd).$$ }

\pf Set $$P_N=\F\times \F[N]=\{ \la x,\eps^{\ell}\ra: x\in \F,\ \ell=0,\dots,N-1\}$$ 
and define the maps 
$$\p(\la x,\eps^{\ell}\ra):=x,\
\f(\la x,\eps^{\ell}\ra):=\la a+x,\eps^{\ell}\ra),\  
\g(\la x,\eps^{\ell}\ra):=\la b+x,\eps^{\ell}\bd(x)\ra.$$ 
One checks easily that the action of $\tilde G$ is well-defined
 and 
that (\ref{p12}) holds.\qed

{\bf Remark} One can easily define in $(\F,+,\cdot, \bd)$ 
functions $\xx,\yy$ and $\zz$ satisfying (\ref{y}) and (\ref{z}).\\

 Assuming that $\F=\C$ and for simplicity that $a\in i\R$ and $b\in \R,$ 
both nonzero, 
we may define, for $z\in \C,$  
$$\bd(z):= \exp(\frac{2\pi i}{N}[{\rm Re}(\frac{z}{a})]).$$ 
This satisfies (\ref{sgn10}) and (\ref{sgn20}) and so $\PPP_N$ over $\C$ is 
definable in $\C$ equipped with the measurable but not continuous 
function above.  \\

{\bf Question} Does there exist a supersimple structure  
 of the form 
 $(\F,+,\cdot, \bd)$ satisfying (\ref{sgn10}) and (\ref{sgn20})?

\epk

\bpk \lb{space2} {\bf The space of semi-definable  
functions}.   
 
Let $\HH$  be an $\F$-algebra containing all the functions
$P_N\to \F$ which are definable in $\PPP_N$ expanded  by  $x,y,z.$\\ 
 
We  define linear operators $\X,\Y,\ZZ,$  
$\FF$ and $\GGG$ on $\HH:$ 
  
\be \lb{*2} \begin{array}{lllll}  
\X: \psi(t)\mapsto \p(t)\cdot \psi(t),\\  
\Y: \psi(t)\mapsto \yy(t)\cdot \psi(t),\\ 
\ZZ: \psi(t)\mapsto \zz(t)\cdot \psi(t),\\ 
\FF:\ \psi(t) \mapsto \psi(\f t),\\  
\GGG:\ \psi(t)\mapsto \psi(\g t).  
\end{array}\ee  
 Denote $\tilde G^*$ the group generated by the  
operators ${\FF},$ ${\GGG}, {\FF}\inv,$ ${\GGG}\inv,$ denote 
$\XX_{\eps}$ (or simply $\XX$) the $\F$-algebra $\F[\X,\Y,\ZZ]$ and $\Ac_{\eps}$ (or simply $\Ac$) the extension of the 
$\F$ algebra $\XX_{\eps}$ by $\tilde G^*.$

While elements of $\HH$ and $\HH$ as a whole are not uniquely defined we prove
in \ref{st2} that 
$\Ac$ is exactly the algebra of operators on $\HH$ generated by $\X,\Y,\ZZ,\FF$ and $\GGG$ satisfying the following defining relations ( ${\bf E}$ stands below for the commutator  
$[\FF,\GGG]$):   
\be \lb{Ac2}\begin{array}{lllllll} 
\X\Y=\Y\X; \X\ZZ=\ZZ\X; \Y\ZZ=\ZZ\Y; \\  
\Y^N=1; \ZZ^N=1; \\  
\FF\X-\X\FF=a\FF; 
\GGG\X-\X\GGG=b\GGG;\\ 
\FF\Y=\eps\Y\FF; 
\Y\GGG=\GGG\Y;\\ 
\ZZ\FF=\FF\ZZ;\\ 
\GGG\ZZ=\Y\ZZ\GGG;\\ 
\FF {\bf E}={\bf E}\FF;   
\GGG {\bf E}= {\bf E} \GGG; 
{\bf E}^N=1.  
\end{array}  
\ee

\epk

\bpk \lb{Irr2} 
Let  
${\rm Max}(\XX)$ be the set of isomorphism classes of  $1$-dimensional  
irreducible $\XX$-modules.\\ 
  
{\bf Lemma 1} {\em ${\rm Max}(\XX)$ can be represented by $1$-dimensional  
 modules $ \la e_{\mu,\xi,\zeta}\ra$ ($ e_{\mu,\xi,\zeta}$ generating the module) 
for $\mu\in \F, \xi,\zeta\in \F[N],$  
defined by the action on the generating vector as follows: 
 $$\X e_{\mu,\xi,\zeta}=\mu e_{\mu,\xi,\zeta},\ \Y 
e_{\mu,\xi,\zeta}=\xi e_{\mu,\xi,\zeta},\ \ZZ e_{\mu,\xi,\zeta}=\zeta e_{\mu,\xi,\zeta}.$$} 
\pf This is a standard fact of commutative algebra.\qed 

{\bf Remark } We can find some of the $e_{\mu,\xi,\zeta}$ in  $\HH,$ 
 which by definition contains the following {\em Dirac delta-functions,} for 
any $p\in P_N,$   
$$\delta_p(t)=\left\{\begin{array}{ll} 1,\ \mbox{ if } t=p;\\ 
                                      0, \mbox{ otherwise}\end{array}\right.$$ 
One checks that 
$$\X\delta_p=\p(p)\cdot\delta_p,\ \ \Y\delta_p=\yy(p)\cdot\delta_p,\ \ 
\ZZ\delta_p=\zz(p)\cdot\delta_p.$$ 
That is we get  $e_{\p(p),\yy(p),\zz(p)}$ in this way. 
\\

Assuming $\F$ is endowed with a fixed function $\bd:  
\F \to \F[N]$   
 we call $\la\mu,\xi,\zeta\ra$ as above {\bf real oriented} if   
$$\bd \mu=\xi.$$  
Correspondingly, we 
call the module $\la e_{\mu,\xi,\zeta}\ra$ real  oriented  
 if $\la\mu,\xi,\zeta\ra$ is. 
 
  ${\rm Max}^{+}(\XX)$ will denote  the   
subspace of ${\rm Max}(\XX)$ consisting of   
real oriented modules $\la e_{\mu,\xi,\zeta}\ra.$\\  
  
{\bf Lemma 2} {\em $\la\mu,\xi,\zeta\ra$ is real oriented if and only if  
 $$\la \mu,\xi,\zeta\ra=\la \p(t),\yy(t),\zz(t)\ra ,$$ 
for some $t\in P_N.$}  \\ \\ 
\pf  It follows from the definition of $\bd$ that $\la \p(t),\yy(t),\zz(t)\ra$ is real oriented.  
 
Assume now that $\la\mu,\xi,\zeta\ra$ is real oriented.  
Since $\p$ is a surjection, there is $t'\in P_N$ such that  
 $\p(t')=\mu.$  
By the definition of $\bd,$ $\yy(t')=\bd\, \mu.$ By the Remark in  
\ref{local2} both values stay the same if we replace $t'$ by $t=[\f,\g]^k\,t'.$ 
By the same Remark, for some $k,$ $\zz(t)=\zeta.$  \qed

Now we introduce an infinite-dimensional $\Ac$-module $\HH_0.$ 
As a vector space $\HH_0$ is spanned by $\{ e_{\mu,\xi,\zeta}:  \mu\in \F, \xi,\zeta\in \F[N]\}.$ The action of the generators of $\Ac$  
on $\HH_0$ is defined on $e_{\mu,\xi,\zeta}$ in accordance with the 
defining relations of $\Ac.$ So, since 
$$\X\FF e_{\mu,\xi,\zeta}=(\FF \X -a\FF)e_{\mu,\xi,\zeta}=  
(\mu-a)\FF e_{\mu,\xi,\zeta},$$ 
$$\Y\FF e_{\mu,\xi,\zeta}=\eps\inv \FF \Y e_{\mu,\xi,\zeta}=  
\eps\inv \xi \FF e_{\mu,\xi,\zeta},$$ 
$$\ZZ\FF e_{\mu,\xi,\zeta}=\FF \ZZ e_{\mu,\xi,\zeta}=  
\zeta\FF e_{\mu,\xi,\zeta},$$ 
and 
$$\X\GGG e_{\mu,\xi,\zeta}=(\GGG \X -b\GGG)e_{\mu,\xi,\zeta}=  
(\mu-b)\GGG e_{\mu,\xi,\zeta},$$ 
$$\Y\GGG e_{\mu,\xi,\zeta}= \GGG \Y e_{\mu,\xi,\zeta}=  
\xi \FF e_{\mu,\xi,\zeta},$$ 
$$\ZZ\GGG e_{\mu,\xi,\zeta}=\GGG\Y\inv \ZZ e_{\mu,\xi,\zeta}=  
\xi\inv\zeta\GGG e_{\mu,\xi,\zeta},$$ 
 
we set $$\FF e_{\mu,\xi,\zeta}:= e_{\f\la \mu,\xi,\zeta\ra},\mbox{ with }  
\f\la \mu, \xi,\zeta\ra=\la\mu-a, \eps\inv\xi,\zeta\ra $$  
and  
 $$\GGG e_{\mu,\xi,\zeta}:= e_{\g\la \mu,\xi,\zeta\ra},\mbox{ with }  
\g\la \mu,\xi,\zeta\ra=\la \mu-b,\xi,\xi\inv\zeta\ra .$$ 
 
From now on we  identify ${\rm Max}^{+}(\XX)$ with 
 the family of real oriented $1$-dimensional  
$\XX$-eigenspaces of $\HH_0.$\\  
   
\bt  \lb{dual2}  
  
(i) There is a bijective correspondence  
$\Xi:\ {\rm Max}^{+}(\XX)\to 
 P_N$  
 between the set 
 of real oriented $\XX$-eigensubspaces of $\HH_0$ and $P_N.$

(ii) The action of  $\tilde G^*$ on $\HH_0$    
induces an action on  ${\rm Max}(\XX)$ and   
leaves ${\rm Max}^{+}(\XX)$ setwise invariant.  
 The correspondence $\Xi$ transfers anti-isomorphically the natural action of   
$\tilde G^*$ on ${\rm Max}^{+}(\XX)$   to a natural action of   
$\tilde G$ on $P_N.$  
  
(iii) The  
map $$\p_{\XX}:\la e_{\mu,\xi,\zeta}\ra\mapsto \mu$$   
is a $N$-to-$1$-surjection ${\rm Max}^{+}(\XX)\to   
\F $ such that  
$$\left( {\rm Max}^{+}(\XX), \FF,\GGG, \p_{\XX},\F \right)\cong_{\xi}  
\left( P_N,\f,\g,\p, \F \right).$$  
\et  
\pf (i) Immediate by Lemma 2. 
 
(ii) Indeed, by the definition above   
the action of $\FF$ and $\GGG$ corresponds  
to the action on real oriented $N$-tuples:  
$$\FF: \la \p(t),\yy(t),\zz(t)\ra \mapsto  
\la  \p(t)-a,\eps\inv \yy(t),\zz(t)\ra=\la \p(\f\inv t),\yy(\f\inv t),\zz(\f\inv t)\ra,$$ 
$$\GGG: \la \p(t)-b,\yy(t),\yy(t)\inv \zz(t)\ra \mapsto  
\la \p(\g\inv t),\yy(\g\inv t),\zz(\g\inv t) \ra.$$  
  
(iii) Immediate from (i) and (ii).\qed  
 
\epk

\bpk \lb{st2}
 {\bf $C^*$-representation.}

Our aim here is to introduce a natural $C^*$-algebra structure on $\Ac.$ In 
fact we will do it for an extension $\Ac^{\#}$ of $\Ac.$ recall that a $\C$-algebra 
$\Ac$ with a norm is called a
$C^*$-algebra if there is a map $*: \Ac\to \Ac$ satisfying the following properties:
for all $X, Y \in \Ac:$
  
 $$(X^*)^* = X,$$
    $$(X Y)^* = Y^* X^*,$$
        $$(X + Y)^* = X^* + Y^*,$$
     for every $\lambda \in \C$ and every $X \in \Ac:$
        $$(\lambda X)^* = \overline{\lambda} X^*$$
 and $$\|X^* X \| = \|X\|^2.$$
In the last condition we do not assume that the norm is always finite.\\

 We will assume $\F=\C,$
$a=\frac{2\pi i}{N},$ $b\in \R$ and
 start by extending the space $\HH$ of semi-definable functions with a 
function $\ww: P_N\to \C$ such that
$$ \exp \ww=\yy,\ \ww(\f t)=\frac{2\pi i}{N}+\ww(t),\  \ww(\g t)=\ww(t). $$
We can easily do this by setting as in (\ref{local2})
$$\ww(\f^m\g^n[\f,\g]^{\ell}\cdot s):=\frac{2\pi im}{N}.$$
Now we extend $\Ac$ to $\Ac^{\#}$ by adding the new operator 
$$\W: \psi\mapsto \ww\psi$$
which obviously satisfies 
$$\W\X=\X\W,\ \W\Y=\Y\W,\ \W\ZZ=\ZZ\W.$$
$$\FF\W=\frac{2\pi i}{N}+\W\FF,\ \GGG\W=\W\GGG.$$

We set 
$$\FF^*:=\FF\inv,\ \GGG^*:=\GGG\inv$$
$$\Y^*:=\Y\inv,\ \W^*:=-\W,\ \X^*:=\X-2\W,$$   
implying that  $\FF,\GGG$ and $\Y$ are unitary and 
$i\W$ and $\X-\W$ are formally selfadjoint.\\

{\bf Lemma} There is a representation of $\Ac^{\#}$ in an inner product space
such that  $\FF,\GGG$ and $\Y$ act as unitary and 
$i\W$ and $\X-\W$ as selfadjoint operators.

\pf Let $\HH_R$ be the subspace of the inner product space $\HH_0$ spanned by vectors
$e_{\mu,\xi,\zeta}$ such that \be \lb{mu1}
\mu=x+\frac{2\pi i k}{N},\ \xi=e^{\frac{2\pi i k}{N}}, \zeta=e^{\frac{2\pi i m}{N}},\mbox{ for }x\in \R,\ k,m\in \Z.\ee
One checks that $\HH_R$ is closed under the
 action of $\Ac$ on $\HH_0$ defined in \ref{Irr2}, that is $\HH_R$ is an $\Ac$-submodule. We also define the action by $\W$
$$\W: e_{\mu,\xi,\zeta}\mapsto \frac{2\pi i k}{N} e_{\mu,\xi,\zeta}$$
for $\mu= x+\frac{2\pi i k}{N}.$ This obviously agrees with the defining relations of $\Ac^{\#}.$ So $\HH_R$ is  an $\Ac^{\#}$-submodule of $\HH_0.$

Now  $\FF$ and $\GGG$ are unitary operators on $\HH_R$ since they transform the
orthonormal basis into itself. $\Y$ is unitary since its eigenvectors form the 
orthonormal basis and the corresponding eigenvalues are of absolute value $1.$ 
  $i\W$ and $\X-\W$ are selfadjoint since their eigenvalues on the orthonormal basis are the reals $-\frac{2\pi k}{N}$ and $x,$ correspondingly.\qed

{\bf Proposition} {\em The action of the algebras $\Ac$ and $\Ac^{\#}$ on
$\HH_R$ are faithful, that is an operator $T$ 
of the algebra annihilates $\HH_R$ if and only if $T=0.$ }

\pf We will prove the statement for $\Ac.$ The proof for  $\Ac^{\#}$ is
similar.

Using the defining relations (\ref{Ac2}) we can represent
$$T=\sum_{i\in I}c_i\X^{i_1}\Y^{i_2}\ZZ^{i_3}\FF^{i_4}\GGG^{i_5}{\bf E}^{i_6}$$
for some finite $I\subset \Z^6,$ $i=\la i_1\ldots i_6\ra$  and 
$c_i\in \C.$

Given an element $e_{\mu,\xi,\zeta}$ of the basis, the action of $T$ on 
it produces
$$Te_{\mu,\xi,\zeta}=\sum_{i\in I}c_i\X^{i_1}\Y^{i_2}\ZZ^{i_3}
e_{\mu(i),\xi(i),\zeta(i)}$$
where  $$e_{\mu(i),\xi(i),\zeta(i)}= \FF^{i_4}\GGG^{i_5}{\bf E}^{i_6}e_{\mu,\xi,\zeta}$$
is a basis element by definition of the action of $\FF$ and $\GGG,$ moreover
one can check that $e_{\mu(i),\xi(i),\zeta(i)}$ are distinct
for distinct $\FF^{i_4}\GGG^{i_5}{\bf E}^{i_6}.$ 

Since the basis elements are eigenvectors of $\X,$ $\Y$ and $\ZZ$  
$$Te_{\mu,\xi,\zeta}=\sum_{i\in I}c_i\cdot d_i(\mu,\xi,\zeta)
e_{\mu(i),\xi(i),\zeta(i)}$$
for some nonzero $d_i(\mu,\xi,\zeta)\in \C.$

Now assume that $T$ annihilates $\HH_R.$ Then the right-hand-side of the above
must be zero and by linear independence all
$c_i\cdot d_i(\mu,\xi,\zeta)=0,$ which can only happen if all $c_i=0$ and $T=0.$\qed  
{\bf Corollary} {\em The $*$-operation on the generators of   $\Ac^{\#}$
defined above extends uniquely to $*$-operation on the whole  $\Ac^{\#}$
and $(\Ac^{\#},\ *)$ satisfies all the identities of a $\C$-algebra with adjoints. Moreover, since $\Ac^{\#}$ has a faithful representation on 
an inner product space we can introduce the usual operator norm on $\Ac^{\#}$
with $\Y,\ZZ,\W,\FF$ and $\GGG$ bounded operators and $\X$ unbounded.  }\\
 
{\bf Remark 1} Our choice of the $C^*$-structure on $\Ac^{\#}$ has been motivated
by 

(i) the need to encode the fact that the relevant $e_{\mu,\xi,\zeta}$ must be
'real oriented', that is $\bd\,\mu=\xi;$

(ii) the natural interpretation of the band function and the related function
$\ww$ (for $a\in i\R$ and $b\in \R$ and $N\to \infty$) as functions indicating
when $\mu$ is 'almost real'. More precisely, as remarked in \ref{angdef2}
$\bd$ can be interpreted, for $a=\frac{2\pi i}{N},$ $b\in \R,$ as 
$$\bd(x+2\pi iy)= \exp 2\pi i \frac{[yN]}{N},$$
where $x,y\in \R,$ and $[yN]$ is the entire part of $yN.$ Since $[yN]/N$
converges to $y$ the condition $\bd\, \mu=1$ says that $\mu$ is 'almost real'.\\

{\bf Remark 2} The natural interpretation of the band function is used in Section~\ref{Hs} to obtain 'the classical limit' $\PPP_{\infty}$ of the $\PPP_N.$ \\

{\bf Comments} 1.We have seen that in the representation $\HH_R$
 the $e_{\mu,\xi,\zeta}$ are eigenvectors of the 
self-adjoint operator $\X-\W.$ So in physics jargon   
$\la  e_{\mu,\xi,\zeta}\ra$
would be called {\bf states}.\\

2. The discrete nature of the imaginary part of $\mu$ in  (\ref{mu1})
is necessitated by two conditions:  the interpretation of $*$ as taking
 adjoints and the noncontinuous form of the band function. The first condition
is crucial for any physical interpretation and the second one follows from the
description of the Zariski structure $\PPP_N.$ Comparing this to the 
real differentiable structure  $\PPP_{\infty}$ constructed in Section~\ref{Hs}
as the limit of the  $\PPP_N$
we suggest to interpret the latter along with its representation via $\Ac$
in this section as the {\em quantisation} of the former.

\epk
\subsection{\bf The algebraic torus case} \lb{intrT}

\bpk    
Let $\F$ be an algebraically closed field of any characteristic prime 
to $N$ an 
$a,b\in \F^*$ be multiplicatively independent.   
 
$G$ acts on $\F^*:$  
$$u x=ax,\ v x=bx.$$  
  
Taking $M$ to be $\F^*$ this determines,  
by subsection~\ref{znorv}, a presmooth  non-algebraic Zariski curve  
$\tilde M$ which from now on we denote $\T_N.$\\

The correspondent definition for the covering map $\p: \tilde M=T_N\to 
 M=\F^*$ 
 then gives us 
\be \lb{p11} \p(\f t)=a\p(t),\ \ \p(\g t)=b\p(t).\ee  
 
We also note that there exists the well-defined function $\p': T_N\to \F$ 
 given  by  
\be \lb{p2} \p'(t)\p(t)=1.\ee 
\epk 
 
For the rest of the section fix $\al=a^{\frac{1}{N}}$ and $\beta=b^{\frac{1}{N}},$ 
roots of $a$ and $b$ of order $N.$ 
 
\bpk \lb{local1} {\bf Semi-definable functions in $\T_N.$}\\  
{\bf Lemma} {\em  
There exist functions    
$$\xx,,\xx',\yy:{T_N}\to \F   $$ 
satisfying  the following {\bf functional equations}, 
 for any $t\in T_N,$  
 
\be \lb{x11}\xx^N(t)=\p(t),\ \xx(\f t)=\al \xx(t),\  \xx(\g t)=
\beta\yy(t)\xx(t)\ee 
\be \xx(t)\xx'(t)=1 \ee 
\be \lb{y11}\yy^N(t)=1,\ \yy(\f t)=\eps \yy(t),\  \yy(\g t)=\yy(t)\ee  
 
} 
 
\pf  Choose a subset $S\subs \F^*$ of representatives of $G$-orbits, 
that is $\F=G\cdot S.$ By construction \ref{tilde}  
we can identify $T_N=\tilde M$ with 
$\tilde G\times S$ in such a way that $\p(\gamma s)=\pr(\gamma)\cdot s.$ 
This means that,  
for any $s\in S$ and $t\in \tilde G\cdot s$ of the form 
$t=\f^m\g^n[\f,\g]^{\ell}\cdot s,$ 
  
$$\p(\f^m\g^n[\f,\g]^{\ell}\cdot s):= 
a^n\cdot b^m\cdot s.$$ 
Fix (randomly) a root $s^{\frac{1}{N}}$ of order $N$ for each $s\in S.$ 
Set $$\xx(\f^m\g^n[\f,\g]^{\ell}\cdot s):= 
\al\cdot \beta\cdot\eps^{-\ell} s^{\frac{1}{N}}.$$  
Set also  
$$\yy(\f^m\g^n[\f,\g]^{\ell}\cdot s):=\eps^{m}.$$ 
 
This satisfies (\ref{x11})- (\ref{y11}). 
\qed  
 
{\bf Remark.} Notice, that it follows from (\ref{x11}) and (\ref{y11}) that 
 
$\xx([\f,\g]t)=\eps\inv \xx(t)$ 
 
$\yy([\f,\g]t)=\yy(t).$

\epk 
 
\bpk \lb{angdef1}  
 
Define the {\bf angular function} on $\F $  
  as a function 
$\sgn: \F  \to \F[N].$ 
  
Set for $\lambda\in \F$  
$$\sgn(\lambda)= \yy(t),\mbox{ if }\p(t)=\lambda.$$ 
This is well-defined by the remark in \ref{local1}.

  
Acting by $\f$ on $t$ and using  (\ref{p11}) and (\ref{y11}) we have 
\be \lb{sgn1} \sgn(a\lambda)=\eps \sgn\, \lambda.\ee  
Acting by $\g$ we obtain 
\be \lb{sgn2}\sgn( b\lambda)=\sgn\, \lambda.\ee

{\bf Proposition}{\em  The structure $\T_N$ is definable in  
$$(\F,+,\cdot, \sgn).$$}
 
Indeed, set $T_N=\F^*$  
and define the maps 
$$\p(t):=t^N$$ 
and  
$$\f(t):=\al t,\  
\g(t):=\sgn(t^N)\beta t.$$ 
One checks easily that $$\g\f( t)=\eps\cdot \f\g( t)$$ and so 
the action of $\tilde G$ is well-defined and 
that (\ref{p11}) holds.\\ 
 
{\bf Remark 1} 
 Assuming that $\F=\C$ and $\eps=\exp(\frac{2\pi i}{N}),$ 
 let for  
an $r\in \R,$ 
 $$a=\exp(\frac{2\pi i}{N}+r),\ \ 
\mbox{ and } 
b\in \R_+,\ b\neq 1.$$ 
Then we may define, for $z\in \C,$  
$$\sgn\,z:= \exp(\frac{2\pi i}{N}[\frac{N}{2\pi}\arg z]).$$ 
This is a well-defined function satisfying also (\ref{sgn1}) and  
(\ref{sgn2}), and so $\T_N$ over $\C$ is 
definable in $\C$ equipped with the measurable but not continuous 
function above. 
 
It is also interesting to remark that, for this angular function, 
$$|\arg z- \frac{2\pi }{N}[\frac{N}{2\pi}\arg z]|\le\frac{2\pi}{N}$$ 
and so $\sgn\, z$ converges uniformly on $z$ to $\exp(i\arg z)$ as $N$ tends to $\infty.$   
 \\

{\bf Remark 2} In  the context of noncommutative geometry 
it is interesting to see whether there exists an abstract, model-theoretic 
interpretation of $\sgn$ which allows a {\em measure theory} for the  
semi-definable functions introduced above. David Evans proved the following  
theorem.\\ 
   
{\bf Theorem} (D.Evans [E]) {\em The class of fields  $(\F,+,\cdot, \sgn)$ 
of a fixed characteristic endowed with a function $\sgn$ 
 satisfying (\ref{sgn1}) and (\ref{sgn2}) has a model companion, which is 
a supersimple  theory. The models of the theory allow a nontrivial finite 
 measure such that all definable sets are measurable.}

\epk

\bpk \lb{space2-1} {\bf The space of semi-definable functions 
and the operator algebra}.   
 
Let $\HH$  be an $\F$-algebra containing all the functions
$T_N\to \F$ which are definable in $\T_N$ expanded  by  $\xx$ and $\yy.$\\ 
 
We  define linear operators $\X,\X\inv,\Y,$  
$\FF$ and $\GGG$ on $\HH:$ 
  
\be \lb{*2-1} \begin{array}{lllll}  
\X: \psi(t)\mapsto \xx(t)\cdot \psi(t),\\ 
\X\inv: \psi(t)\mapsto \xx'(t)\cdot \psi(t),\\ 
\Y: \psi(t)\mapsto \yy(t)\cdot \psi(t),\\ 
\FF:\ \psi(t) \mapsto \psi(\f t),\\  
\GGG:\ \psi(t)\mapsto \psi(\g t).  
\end{array}\ee 
 Denote $\tilde G^*$ the group generated by the  
operators ${\FF},$ ${\GGG}, {\FF}\inv,$ ${\GGG}\inv,$ denote 
$\XX_{\eps}$ the $\F$-algebra $\F[\X,\X\inv,\Y]$ and  
$\Ac_{\eps}$ (or simply $\Ac$) the extension of the 
$\F$ algebra $\XX_{\eps}$ by $\tilde G^*.$

The generators of the algebra $\Ac_{\eps}$ obviously satisfy 
the following relations, for  ${\bf E}$ standing for the commutator  
$[\FF,\GGG],$   
\be \lb{Ac1} \begin{array}{lllllllll} 
\X\Y=\Y\X; \\  
\Y^N=1;\ \X\X\inv=1;\\  
\X\FF=\al\inv\FF\X;\\ 
\X\GGG=\beta\inv\Y\inv \GGG\X;\\ 
\Y\FF=\eps\inv\FF\Y;\\ 
\Y\GGG=\GGG\Y;\\ 
\FF {\bf E}={\bf E}\FF;   
\GGG {\bf E}= {\bf E} \GGG; 
{\bf E}^N=1.  
\end{array}  
\ee  
 
We prove later on, 
in the Proposition and Corollary of \ref{st2}, that the algebra
determined by the relations (\ref{Ac1})
is exactly $\Ac$ and so the definition of $\Ac$ does not depend on the arbitrariness in the construction of $\HH.$  

\epk  
 
\bpk \lb{Irr1} 
Let  
${\rm Max}(\XX)$ be the set of isomorphism classes of  $1$-dimensional  
irreducible $\XX$-modules.\\ 
  
{\bf Lemma 1} {\em ${\rm Max}(\XX)$ can be represented by $1$-dimensional  
 modules $ \la e_{\mu,\xi}\ra$ ($=\F e_{\mu,\xi}$) 
for $\mu\in \F, \xi\in \F[N],$  
defined by the action on the corresponding generating vector: 
 $$\X e_{\mu,\xi}=\mu e_{\mu,\xi},\ \Y 
e_{\mu,\xi}=\xi e_{\mu,\xi}.$$} 
\pf This is a standard fact of commutative algebra.\qed

Now we introduce an infinite-dimensional $\Ac$-module $\HH_0.$ 
As a vector space $\HH_0$ is spanned by $\{ e_{\mu,\xi}:  
 \mu\in \F, \xi\in \F[N]\}.$ The action of the generators of $\Ac$  
on $\HH_0$ is defined on $e_{\mu,\xi}$ in accordance with the 
defining relations of $\Ac.$ So, since 
$$\X\FF e_{\mu,\xi}=\al\inv\FF \X e_{\mu,\xi}=  
\al\inv\mu\FF e_{\mu,\xi},$$ 
$$\Y\FF e_{\mu,\xi}=\eps\inv\FF \Y e_{\mu,\xi}=  
\eps\inv\xi \FF e_{\mu,\xi},$$ 
and 
$$\X\GGG e_{\mu,\xi}=\beta\inv\GGG\Y\inv \X e_{\mu,\xi}= 
\beta\inv\xi\inv\mu\GGG e_{\mu,\xi},$$ 
$$\Y\GGG e_{\mu,\xi}=\GGG \Y e_{\mu,\xi}= 
\xi\GGG e_{\mu,\xi},$$ 
we set $$\FF e_{\mu,\xi}:= e_{\nu,\zeta},\mbox{ with }  
\nu=\al \mu,\ \ \zeta=\eps\inv\xi$$  
and 
 $$\GGG e_{\mu,\xi}:= e_{\nu,\zeta},\mbox{ with }  
\nu=\beta\xi\inv\mu,\ \ \zeta=\xi .$$  
 
We may now identify ${\rm Max}(\XX)$ as the family of $1$-dimensional   
$\XX$-eigenspaces of $\HH_0.$ \\

Assuming $\F$ is endowed with an angular function $\sgn$ 
 we call $\la\mu,\xi\ra$ as above {\bf positively oriented} if   
$$\sgn\, \mu^N=\xi.$$  
Correspondingly, we 
call the $\XX$-module (state) $\la e_{\mu,\xi}\ra$ positively oriented  
 if $\la\mu,\xi\ra$ is. 
  $\HH_0^{+}$ will denote  the linear subspace of   
 $\HH_0$ spanned by the  
positively oriented states $\la e_{\mu,\xi}\ra.$
We denote  ${\rm Max}^+(\XX)$ the family of $1$-dimensional positively 
oriented  
$\XX$-eigenspaces of $\HH_0,$ or {\em states} as such things are
referred to in 
physics literature.
\\ \\ 
{\bf Lemma 2} {\em (i) $\la\mu,\xi\ra$ is positively oriented if and only if  
 $$\la \mu,\xi\ra=\la \xx(t),\yy(t)\ra ,$$ 
for some $t\in T.$ 
 
(ii) $\HH_0^{+}$ is invariant under the action of $\FF$ and  
$\GGG,$ so is an $\Ac$-module. }
 \\ \\ 
\pf (i) 
Indeed, since $p$ is a surjection, there is $t'\in T$ such that  
 $\p(t')=\mu^N.$ Hence, by the definition of $\xx(t')$ and $\sgn(t')$  
we have $\xx(t')=\eps^k\mu,$ $\yy(t')=\xi,$ 
for some $k.$ By the Remark in \ref{local1} we have $\xx([\f,\g]^kt')=\eps^{-k}\xx(t')=\mu$ and  $\yy([\f,\g]^kt')=\yy(t')=\xi.$ So  $t=[\f,\g]^kt'$ is as required.

(ii)  Immediate by the definition of the action.\qed 
{\bf Remark 2} It is immediate from the Lemma and  
Remark 1 that all the positively oriented 
$e_{\mu,\xi}$ are represented by the Dirac functions $\delta_t,$ $t\in T_N.$

\bt  \lb{dual1} 
  
(i) There is a bijective correspondence  
$\Xi:\ {\rm Max}^{+}(\XX)\to 
 T_N$  
 between the set 
 of positively oriented states and $T_N.$  
 
(ii) The action of  $\tilde G^*$ on $\HH_0$    
induces an action on  ${\rm Max}(\XX)$ and  
leaves ${\rm Max}^{+}(\XX)$ setwise invariant. 
 The correspondence $\Xi$ transfers anti-isomorphically the natural action of   
$\tilde G^*$ on ${\rm Max}^{+}(\XX)$   to the natural action of   
$\tilde G$ on $T_N.$ 
  
(iii) The  
map $$\p_{\XX}:\la e_{\mu,\xi}\ra\mapsto \mu^N$$  
is a $N$-to-$1$-surjection ${\rm Max}^{+}(\XX)\to   
\F $ such that  
$$\left( {\rm Max}^{+}(\XX), \FF,\GGG, \p_{\XX},\F \right) 
\cong_{\Xi}  
\left( T_N,\f,\g,\p, \F \right).$$  
\et  
\pf (i) Immediate by Lemma 2 of \ref{Irr1}. 
  
(ii) Indeed, by the definition above  
the action of $\FF$ and $\GGG$ corresponds 
to the action on positively oriented pairs: 
$$\f: e_{ \xx(t),\yy(t)} \mapsto \FF\inv e_{ \xx(t),\yy(t)}= 
e_{\al \xx(t),\eps \yy(t)}=e_{ \xx(\f t),\yy(\f t)},$$ 
$$\g: e_{ \xx(t),\yy(t)} \mapsto \GGG\inv e_{ \xx(t),\yy(t)}= 
e_{\beta\yy(t)\xx(t),\ \yy(t)}=e_{ \xx(\g t),\yy(\g t)}.$$ 
 
(iii) Immediate from (i) and (ii).\qed 
 
\epk  
\bpk  {\bf $C^*$-structure}

We add to \ref{local1} the new semi-definable function $\ww$ satisfying,
for some $\delta,$ such that $\delta^N=\eps,$
$$\yy=\ww^N,\ \ \ww(\f t)=\delta \ww(t),\ \ \ww(\g t)=\yy(t)\inv\ww(t)\ \ 
.$$
In accordance with \ref{local1} we can define 
$$\ww(\f^m\g^n[\f,\g]^l)=\delta^m\eps^{l}.$$
Now we introduce
$$\hat\sgn\,x:=\ww(t),\ \mbox{ for }x=\xx(t).$$
Since $\xx$ is a bijection  this is well-defined on $\F.$
Moreover, using the unique representation  
$$x=\xx(\f^m\g^n[\f,\g]^l)=\al^m\beta^n\eps^{-l}s^{\frac{1}{N}}$$
of \ref{local1} we have $$\hat\sgn(\al^m\beta^n\eps^{-l}s^{\frac{1}{N}})= \delta^m\eps^{l}.$$
\\

Taking $a=\eps\rho,$ $\rho\in \R_+$ (positive reals), $\rho\neq 1,$ as suggested in \ref{angdef1},
 and $\alpha\inv\delta\in \R_+,$  we have
$$\hat \sgn(\alpha x)=\delta \hat \sgn\,x,\ \ \ \hat \sgn(\beta x)=
\hat \sgn\,x.$$

Extend the list of operators on $\HH$ to include
$$\W: \psi \mapsto \ww\cdot\psi.$$
Obviously $\W$ commutes with $\X.$ As in \ref{Irr1} denote $e_{\mu,w}$ an eigenvector of $\X$ and $\W$ with eigenvalues $\mu$ and $w$ correspondingly. The action of $\FF$ and $\GGG$ is defined on $e_{\mu,w}$ similarly to \ref{Irr1}:
\be \lb{FG}
\begin{array}{ll} \FF: e_{\mu,w}\mapsto  e_{\alpha\inv\mu,\,\delta\inv w}\\  
 \GGG: e_{\mu,w}\mapsto  e_{\beta\inv w^{-N}\mu, w^{1-N}}\end{array}
 \ee

Consider the algebra $\Ac$ as a $C^*$-algebra with the 
 condition that $\X\W\inv$ is 
{\bf selfadjoint} and $\W,$ $\FF$ and $\GGG$ are {\bf unitary}.

Set
 $$\W^*:=\W\inv,\ \FF^*:=\FF\inv,\ \GGG^*:=\GGG\inv $$
that is define these operators as {\em unitary}.
Set   $$\X^*:=\W\inv\X\W\inv
=\X\W^{-2},$$
so
$$(\X\W\inv)^*={\W^*}\inv\X^*=\W\X^*=\X\W\inv$$
 that is $\X\W\inv$ is {\em selfadjoint}.

Of course $$[\FF,\GGG]^*=[\FF,\GGG]\inv$$ so $[\FF,\GGG]$
 is unitary as well.\\

{\bf Lemma} 
There is an inner product space $\HH_+$ with the faithful action of
   $\Ac$ on it such that $*$ corresponds to taking adjoint operators.

\pf Consider $\HH_+\subs \HH$ generated by all $e_{\mu,w}$ 
satisfying the condition
\be \lb{muw}\mu\cdot w\inv\in \R_+,\ w=\exp \frac{2\pi i k}{N^2},\mbox{ for }k\in \Z.\ee

We introduce the inner product in $\HH_+$
assuming  the $e_{\mu,w}$ to form an orthonormal basis.

Now, by definition $\X\W\inv$ acts as a positive selfadjoint operator
$$ \X\W\inv: e_{\mu,w}\mapsto  \mu w\inv e_{\mu,w}.$$
$\W$ acts as  unitary since $w$ is a root of unity.

 $\HH_+$ is closed under $\FF$ and $\GGG$ since
 $\alpha\inv\mu\delta w\inv$ and  
$ \beta\inv\mu\delta w\inv$ are in $\R_+.$

The fact that the action is faithful (that is the only operator that annihilates $\HH$ is $0$) is essentially proved in  the Proposition and Corollary of \ref{st2}.
\qed

{\bf Comment} Using the representation on $\HH_+$ 
one clearly can interpret the angular function $\hat \sgn\, \mu$ as 
$\exp \arg \mu,$ for $\mu$ satisfying (\ref{muw}). For general $\mu$ 
 we can use  the interpretation  as in \ref{angdef1}:
$$\hat\sgn\,\mu=\exp \frac{2\pi i}{N^2}[\frac{N^2}{2\pi}\arg \mu],$$
where $[r]$ stands for the integer part of a real number $r.$
Of course, we stress again that  $\hat\sgn\,\mu$ is very well approximated by
$\exp \arg \mu$ since
 $$|\frac{2\pi i}{N^2}[\frac{N^2}{2\pi}\arg \mu]-\arg \mu|\le 
\frac{2\pi}{N^2}.$$

In other words, the condition on the states being positively oriented in
Theorem~\ref{dual1} is similar to conditions
usually stated in terms of $C^*$-algebras.
 This must justify the name
$*$-functions for $\sgn,$ $\hat\sgn$ and $\bd.$ 
\epk
\section{The metric limit}
 \lb{Hs}
Our aim in this section is to find an interpretation of the limit, as
$N$ tends to $\infty,$
of structures $\T_N$ or $\PPP_N$ in ``classical'' terms.
``Classical'' here is supposed to mean `` using 
function and relations given in terms of real manifolds and analytic
 functions''. 
Of course, we
have to define the meaning of the ``limit'' first. We found a
satisfactory
solution to this problem in case of $\PPP_N$ which is presented below.
 
\bpk
First we want to establish a connection of the group $\tilde G_N$ 
with  the 
{\em integer Heisenberg group} 
$H(\Z)$ which is the group of matrices of the form
\be\lb{Heis}\left(\begin{array}{lll} 1\ \ \ k\ \ \ m\\
                           0\ \ \ 1\ \ \ l\\
                           0\ \ \ 0\ \ \ 1
\end{array}
\right)\ee
with $k,l,m\in \Z.$ More precisely,  $\tilde G_N$ is isomorphic to
the group $$H(\Z)_N=H(\Z)/ N. Z,$$ where
$N.Z$ is the central subgroup
$$N.Z=\left\{\ \left(\begin{array}{lll} 1\ \ \ 0\ \ Nm\\
                           0\ \ \ 1\ \ \ 0\\
                           0\ \ \ 0\ \ \ 1
\end{array}
\right):\ m\in \Z\right\}
$$ 

Similarly   the 
{\em real Heisenberg group} 
$H(\R)$ is defined as the group of matrices of the form (\ref{Heis})
with $k,l,m\in \R.$ 
The analogue (or the limit case) of $H(\Z)_N$ is the factor-group
 $$H(\R)_{\infty}:=\ H(\R)/
\left(\begin{array}{lll} 1\ \ \ 0\  \ \Z\\
                           0\ \ \ 1\ \ \ 0\\
                           0\ \ \ 0\ \ \ 1
\end{array}
\right)$$

In fact there is the natural group embedding
$$i_N: \left(\begin{array}{lll} 1\ \ \ k\ \ \ m\\
                           0\ \ \ 1\ \ \ \ell\\
                           0\ \ \ 0\ \ \ 1
\end{array}
\right) \mapsto \left(\begin{array}{lll} 1\ \  {\frac{k}{\sqrt{N}}} \ \ \ \frac{m}{N}\\
                           0\ \ \ 1\ \  \ \  {\frac{\ell}{\sqrt{N}}}\\
                           0\ \ \ 0\ \ \ \ 1
\end{array}
\right)$$
inducing the embedding $H(\Z)_N\subset H(\R)_{\infty}.$\\

Notice the following\\

{\bf Lemma 1} {\em Given the embedding $i_N$ for every 
$\la u,v,w\ra\in H(\R)_{\infty}$ there is $\la \frac{k}{\sqrt{N}},
\frac{\ell}{\sqrt{N}},\frac{m}{{N}}\ra\in i_N(H(\Z)_N)$
such that 
$$|u-\frac{k}{\sqrt{N}}|+|v-\frac{\ell}{\sqrt{N}}|+|w-\frac{m}{{N}}|
<\frac{3}{\sqrt{N}}.$$}

In other words, the distance (given by the sum of absolute values)
between any point of $H(\R)_{\infty}$ and the set $i_N(H(\Z)_{N})$ is at
most $3/\sqrt N.$ Obviously, also the distance  between any point of 
$i_N(H(\Z)_{N})$ and the set $H(\R)_{\infty}$ is $0,$ because of the
embedding. In other words, this defines that {\bf the Hausdorff distance between the
two sets is at most  $3/\sqrt N.$}

In situations when the pointwise distance between sets $M_1$ and $M_2$
is defined we also say that the Hausdorff distance between two 
$L$-structures on
$M_1$ and $M_2$ is at most $\alpha$ if the
 Hausdorff distance between the universes $M_1$ and $M_2$ as
well as between $R(M_1)$ and $R(M_2),$ for any $L$-predicate  or
graph of an $L$-operation $R,$ is
at most $\alpha.$

Finally, we say that an $L$-structure $M$ is the {\bf  Hausdorff
limit}
of $L$-structures $M_N,$ $N\in \N,$ if for each positive $\alpha$
there is $N_0$ such that for all $N>N_0$ the distance between $M_N$
and
$M$ is at most $\alpha.$  \\ \\
{\bf Remark} It makes sense to consider the similar notion of Gromov-Hausdorff
distance and Gromov-Hausdorff limit.\\ \\
{\bf Lemma 2} {\em The group structure  $H(\R)_{\infty}$ is the
Hausdorff
limit of its substructures $H(\Z)_{N},$ where the distance is defined
by
the embeddings $i_N.$}

\pf Lemma 1 proves that the universe of  $H(\R)_{\infty}$ is the limit
of the corresponding sequence. Since the group operation is continuous
in the topology determined by the distance, the graphs of the group
operations converge as well.\qed

\epk

\bpk \lb{Hac}

Given nonzero real numbers $a,b,c$  
the integer Heisenberg group $H(\Z)$ acts on $\R^3$ as follows:
\be \lb{Hact}
\la k,l,m\ra \la x,y,s\ra=\la x+ak,y+bl,s+acky+abcm\ra
\ee 
where $\la k,l,m\ra$ is the matrix (\ref{Heis}).

We can also define the action of $H(\Z)$ on $\C\times S^1,$
equivalently on $\R\times \R\times \R/\Z,$ as follows
\be \lb{Hact1}
\la k,l,m\ra\la x,y,\exp {2\pi is}\ra=\la x+ak,y+bl,\exp 2\pi i(s +acky+abcm)\ra
\ee 
where $x,y,s\in \R.$

In the discrete version intended to model \ref{intrT2} we consider
 $\frac{q}{N},$ $q\in \Z,$ in place
of $s\in \R$ and take $a=b=\frac{1}{\sqrt N}.$  We replace
 (\ref{Hact1}) by

\be \lb{Hact1d0}
\la k,l,m\ra\la x,y,e^{ \frac{2\pi iq}{N}}\ra=
\la x+\frac{k}{\sqrt{N}},\ y+\frac{\ell}{\sqrt N},\ \exp {2\pi i\frac{q+k[y\sqrt N]+m}{N}}\ra
\ee

Check that this is still an action:
$$\la k',\ell',m'\ra(\la k,\ell,m\ra\la x,y,\exp \frac{2\pi iq}{N}\ra)=$$
$$\la k',\ell',m'\ra\la x+\frac{k}{\sqrt{N}},y+\frac{\ell}{\sqrt N},\exp2\pi i\frac{q+k[y\sqrt N]+m}{N}\ra=$$
$$\la x+\frac{k}{\sqrt{N}}+\frac{k'}{\sqrt N},y+\frac{\ell}{\sqrt
N}+\frac{\ell'}{\sqrt N},
\exp2\pi i\frac{q+k[y\sqrt N]+m+k'[(y+\frac{\ell}{\sqrt N})\sqrt N]+m'}{N}\ra=$$
$$\la x+\frac{k+k'}{\sqrt N},y+\frac{\ell+\ell'}{\sqrt N},
\exp2\pi i\frac{q+(k+k')[y\sqrt N]+k'l+m+m'}{N}\ra=$$
$$(\la k',l',m'\ra\la k,l,m\ra)\la x,y,\exp \frac{2\pi iq}{N}\ra$$

Moreover, we may take $m$ modulo $N$ in (\ref{Hact1d0}), that is 
$\la k,l,m\ra\in H(\Z)_N,$ and
 simple calculations similar to the above show the following\\

{\bf Lemma 1} {\em The formula (\ref{Hact1d0}) defines
the free action of $H(\Z)_N$ on $\R\times \R\times \exp \frac{2\pi i}{N}
\Z$ (equivalently on 
$\C\times \exp \frac{2\pi i}{N}\Z$)}.\\

We think of $\la x,y,\exp \frac{2\pi iq}{N}\ra$ as an element $t$
of $P_N$ (see \ref{intrT2}), $x+iy$ as $p(t)\in \C.$
The actions $x+iy\mapsto a+x+iy$ and  $x+iy\mapsto x+i(y+b)$
are obvious rational automorphisms of the affine line $\C.$

We interpret the action of $\la 1,0,0\ra$ and $\la 0,1,0\ra$
by  (\ref{Hact1d0}) on $\C\times \exp \frac{2\pi i}{N}\Z$
as $\f$ and $\g$ correspondingly. Then the commutator $[\f,\g]$
corresponds to $\la 0,0,-1\ra,$ which is the generating element of
the centre of $H(\Z)_N.$ In other words, the subgroup ${\rm gp}( \f,\g)$
of $H(\Z)_N$ generated by the two elements is isomorphic to $\tilde G.$ 
We thus get, using Lemma 1 of \ref{tilde}\\

{\bf Lemma 2} {\em Under the above assumption and notation
the structure on $\C\times  \exp \frac{2\pi i}{N} \Z$ in the language
of \ref{intrT2}
described by (\ref{Hact1d0}) 
is isomorphic to the example $\PPP_N$ of
 \ref{intrT2} with $\F=\C.$}\\

Below we identify $\PPP_N$ with the structure above based on $\C\times
\{ \exp \frac{2\pi i}{N}\Z\}.$ 

Note that every group word in $\f$ and $\g$ gives rise to a definable
map in $\PPP_N.$ We want introduce a uniform notation for such
definable functions.
 \\

Let $\al$ be a monotone nondecreasing converging sequence of the form
$$\al=\{\frac{k_N}{\sqrt N}:  k_N, N\in \Z,\ N>0\}.$$
We call such a sequence {\bf admissible} if there is an $r\in \R$ such that
\be \lb{k/} |r-\frac{k_N}{\sqrt N}|\le \frac{1}{\sqrt N}.\ee  
Given $r\in \R$ and $N\in \N$ one can easily find $k_N$ satisfying 
(\ref{k/}) and so construct an $\al$ converging to $r,$ which we denote $\hat \al,$
$$\hat \al:=\lim \al=\lim_N \frac{k_N}{\sqrt N}.$$ 

We denote $I$ the set of all admissible sequences converging to a real
on $[0,1],$ so
$$\{ \hat\al:\al\in I\}=\R\cap [0,1].$$

For each $\al\in I$ we introduce two operation symbols $\f_{\al}$ and
$\g_{\al}.$ We denote $\PPP^{\#}_N$ the definable
 expansion of $\PPP_N$ by all such
symbols with the interpretation
 
$$\f_{\al}=\f^{k_N},\ \ \ \g_{\al}=\g^{k_N}\mbox{\ \ ($k_N$-multiple of the operation)},$$
if $\frac{k_N}{\sqrt N}$ stands in the $N$th position in
 the sequence $\al.$

Note that the sequence 
$$dt:=\{\frac{1}{\sqrt N}:  N\in \N\}$$
is in $I$ and $\f_{dt}=\f,$ $\g_{dt}=\g$ in all $\PPP^{\#}_N.$

\epk

\bpk

We now define the structure $\PPP_{\infty}$ 
to be the structure on sorts
$\C\times S^1$ (denoted $P_{\infty}$) and sort 
$\C,$ with the field structure on $\C$ and the projection map
$\p: \la x,y,e^{2\pi is}\ra \mapsto \la x,y\ra\in \C,$
and  definable maps $\f_{\al}$ and $\g_{\beta},$ $\al,\beta\in I,$ acting
on $\C\times S^1$ (in accordance with the action by $H(\R)_{\infty}$)
as follows 
\be \lb{Hactinf}\begin{array}{ll}
\f_{\al}(\la x,y,e^{2\pi is}\ra)=\la \hat \al,0,0\ra\la x,y,e^{2\pi is}\ra=
\la x+\hat \al,\ y,\
e^{2\pi i(s +\hat \al y)}\ra\\

\g_{\beta}(\la x,y,e^{2\pi is}\ra)=\la 0,\hat \beta,0\ra\la x,y,e^{2\pi is}\ra=
\la x,\ y+\hat \beta,\ 
e^{2\pi is }\ra\end{array}
\ee

{\bf Theorem 1} {\em $\PPP_{\infty}$ is the Hausdorff limit
of structures $\PPP^{\#}_N.$ }

\pf  The sort $\C$ is the same in all structures.

The sort $P_{\infty}$ is the limit of its substructures $P_N$ since $S^1$ ($=\exp i\R$) is the limit of
$ \exp \frac{2\pi i}{N}\Z$ in the standard metric of $\C.$ Also, the
graph of
the projection map $\p:\ P_{\infty}\to \C$  is the limit of  $\p:\
P_{N}\to \C$
for the same reason.   

Finally it remains to check that the graphs of $\f$ and $\g$ in
$\PPP_{\infty}$  are the limits of those in $\PPP_N.$ It is enough
to see that for any 
 $\la x,y,\exp \frac{2\pi iq}{N}\ra\in
P_N$ the result of the action by $\f_{\al}$ and $\g_{\beta}$
calculated in 
$ \PPP^{\#}_N$ is at most at the distance $2/\sqrt N$ from
the ones calculated in $\PPP_{\infty},$ for any 
 $\la x,y,\exp \frac{2\pi iq}{N}\ra\in P_{\infty}.$
And indeed, the action in
$\PPP_N^{\#}$ by definition is  
\be \lb{act}\begin{array}{ll}\f_{\al}: 
\la x,y,\exp \frac{2\pi iq}{N}\ra\mapsto 
\la x+\frac{k_N}{\sqrt{N}},\ y,\ \exp\frac{2\pi i}{N}(q+k_N[y\sqrt N])\ra\\
\g_{\beta}: \la x,y,\exp \frac{2\pi iq}{N}\ra\mapsto  
\la x,\ y+\frac{\l_N}{\sqrt{N}},\ \exp{2\pi i\frac{q}{N}}\ra

\end{array}\ee


Obviously,
$$|\frac{k_Ny}{\sqrt N}-\frac{k_N[y\sqrt N]}{N}|=\frac{k_N}{\sqrt
N}|\frac{y\sqrt N-[y\sqrt N]}{\sqrt N}|< \frac{k_N}{\sqrt
N}\frac{1}{\sqrt N}
\le \frac{1}{\sqrt N},$$
which together with (\ref{k/}) proves that the right hand side of (\ref{act}) is at the
distance at most $\frac{2}{\sqrt N}$ from the right hand side of 
(\ref{Hactinf}) uniformly on the point  $\la x,y, \exp\frac{2\pi i
q}{N}\ra.$
\qed
\epk

\bpk \lb{nabla}  The structure $\PPP_{\infty}$ can be seen as the
principal bundle over $\R\times \R$ with the structure group $U(1)$ (the
rotations of $S^1$) and the projection map $\p.$ The action by 
the Heisenberg group allows to define a {\em connection} on the bundle.
A connection determines ``a smooth transition from a point in a fibre to a point 
in a nearby fibre''.
As noted above  $\f$ and $\g$ in the limit process correspond to 
infinitesimal actions (in a nonstandard model of  $\PPP_{\infty}$)
  which can be written
in the 
form  
$$\begin{array}{ll}
\f(\la x,y,e^{2\pi is}\ra)=
\la x+dx,\ y,\
e^{2\pi i(s +ydx)}\ra\\

\g(\la x,y,e^{2\pi is}\ra)=
\la x,\ y+dy,\ 
e^{2\pi is }\ra\end{array}
$$   
where $dx$ and $dy$ are infinitesimals equal to the $dt$ of \ref{Hac}.

These formulas allow to calculate the derivative of a section
$$\psi: \la x,y\ra\mapsto \la x,y,e^{2\pi is(x,y)}\ra$$
of the bundle in any direction on $\R\times \R.$
In general moving infinitesimally from the point $\la x,y\ra$ along $x$ we get 
$\la x+dx,\ y,\ \exp 2\pi i (s+ds)\ra.$ We need to compare this to
 {\em the parallel transport along $x$} given by the formulas above, 
$\la x+dx,\ y,\ \exp 2\pi i (s+ydx)\ra.$ So the difference is
$$\la 0,\ 0,\ \exp 2\pi i(s+ ds)-\exp 2\pi i(s+ ydx)\ra.$$
Using the usual laws of differentiation one gets for the third term
$$\begin{array}{lll}\exp 2\pi i(s+ ds)-\exp 2\pi i(s+ ydx)=\\
(\exp 2\pi i(s+ ds)-\exp 2\pi is) -(\exp 2\pi i(s+ ydx)-\exp 2\pi is)=\\
d \exp 2\pi is -2\pi i y \exp 2\pi is\, dx=
(\frac{d \exp 2\pi is}{dx} -2\pi i y \exp 2\pi is)dx\end{array}$$
which gives for a section $\psi=\exp 2\pi is$ the following {\em covariant derivative} along $x,$
$$\nabla_x \psi=\frac {d }{dx}\psi -2\pi iy\psi=(\frac {d }{dx}+A_x)\psi.$$ 
Similarly, the covariant derivative along $y$ 
$$\nabla_y\psi=\frac {d }{dy}\psi= (\frac {d }{dy}+A_y)\psi$$
 with the second term $A_y=0$

The {\em curvature} of the connection is by definition the commutator
$$[\nabla_x,\nabla_y]=\frac{d A_y}{dx}-\frac{d A_x }{dy}=2\pi i,$$
that is 
in physicists terms this pictures an $U(1)$-gauge field theory over $\R^2$
with a constant nonzero curvature. 
\epk

\subsection{Algebraic torus}

\bpk
We think of elements of  $\C^*\times S^1$ as pairs $\la z, \exp is\ra,$ 
where $z=\exp(ix+y)\in \C^*$ $x,y,s\in \R.$

The action of  $H(\Z)$ on  
 $\C^*\times S^1,$ can be given, following (\ref{Hact}) by

\be \lb{Hact3}\begin{array}{ll}
\f 
\la\exp (ix+y),\ \exp {is} \ra 
=\la\exp (ix+ia+y),\ \exp i(s+ay)\ra
\\
\g 
\la\exp ix+y,\ \exp {is}\ra 
=\la\exp (ix+y+b),\ \exp is\ra
\end{array}
\ee 
 
The action by $\g$ is well-defined since it simply takes the pair $\la z,t\ra$ 
to $\la e^bz,t\ra.$

To calculate $\f\,\la z,t\ra$ one first takes $$\ln z=ix+y+2\pi in=i(x+2\pi n)
+ y,\ \ n\in \Z.$$

This recovers $y$ uniquely and so $\f$ is well-defined.

The corresponding discrete version will be
\be \lb{Hact3d}
\la k,l,m\ra \la\exp (2\pi ix+y),\ \exp{2\pi i \frac{q}{N}}\ra=$$ 
$$=\la\exp (2\pi i(x+\frac{k}{N})+y+\frac{\ell}{N}),\ \exp 2\pi i\frac{q+k[Ny]+m}{N}\ra
\ee 
This is a group action, by the same calculation as in \ref{Hac}.

In this discrete version $t=\la\exp (2\pi ix+y),\ \exp{2\pi i \frac{q}{N}}$ is an element of $T_N$ and correspondingly $\p(t)=\exp (2\pi ix+y).$ 
The $a$ and $b$
of  \ref{intrT} will be $e^{\frac{2\pi i}{N}}$ and 
$e^{\frac{1}{N}}$ correspondingly.

\bt The structure on $\C^*\times \{ \exp \frac{2\pi i \Z}{N}\}$ in the language
of \ref{intrT}
described by (\ref{Hact3d}) 
is isomorphic to the example $\T_N$ of
 \ref{intrT} with $\F=\C.$\et

We want to calculate the covariant derivative following the method of 
\ref{nabla}. We use similar notation for the infinitesimal action 
$$dx=\{ \frac{1}{N}: N\in \N\}=dy,$$
the infinitesimal corresponding to the sequence.
But the actual cooridinates on $\C^*$ are 
$$z^1=e^{2\pi ix}\mbox{ and }z^2=e^y,$$
so $$dz^1=2\pi i z^1dx,\ \ dz^2=z^2dy.$$

Now for $$\psi: z\mapsto \la x,y,e^{2\pi is(z)}\ra$$
the difference between the shift $dz^1$ and the parallel transport along the same shift will be, by the same formulas as in \ref{nabla},
$$\exp 2\pi i(s+ ds)-\exp 2\pi i(s+ ydx).$$
This is equal to
$$\begin{array}{ll}
(\frac{d \exp 2\pi is}{dx} -2\pi i y \exp 2\pi is)\, dx=\\
(\frac{d \exp 2\pi is}{dz^1} -\frac{\ln z^2}{z^1} \exp 2\pi is)\,dz^1
\end{array}$$
which gives the covariant derivative along $z^1$
$$\nabla_{z^1} \psi=\frac {d }{dz^1}\psi -\frac{\ln z^2}{z^1}\psi.$$ 
Similarly, $\nabla_{z^2}$ the covariant derivative along $z^2$ 
  is just
$\frac{d }{dz^2}\psi,$ the second term zero.

The curvature of the connection is
$$[\nabla_{z^1},\nabla_{z^2}]= \frac{1}{z^1z^2},$$
which is  a nonconstant  curvature (note also that $z^1z^2=\exp(2\pi i x+ y)$ does not vanish
on $\C^*$).

\epk

{\bf References}\\

[DFN] B.A. Dubrovin, A.T. Fomenko, S.P. Novikov, {\bf Modern geometry: methods and applications}, New York; London: Springer-Verlag, 1990\\

[HZ]  E.Hrushovski and B.Zilber, {\em Zariski  Geometries}. Journal of AMS, 9(1996),
1-56  \\

[S] G. Svetlichny, {\bf Preparation for Gauge Theory} arXiv:math-ph/9902027v1\\

[Z0] B.Zilber, {\em Lecture notes on Zariski structures}, 1994-2006, web-page\\

[Z1] B.Zilber, {\em A class of quantum Zariski geometries}, to appear in 
Proc. Newton Institute 2005 Programme.

\end{document}